\documentclass[reqno,12pt]{article}

\usepackage{amsfonts}
\usepackage{graphicx}
\usepackage{epstopdf}
\usepackage{amsmath}
\usepackage{amsthm}
\usepackage{latexsym,bm}
\usepackage{amssymb}
\usepackage{indentfirst}
\usepackage{hyperref}

\usepackage{setspace}
\newtheorem{thm}{Theorem}[section]
\newtheorem{cor}[thm]{Corollary}

{\theoremstyle{definition}

\newtheorem*{conjecture}{Conjecture}
}

\numberwithin{equation}{section}

\begin{document}
\title{\bf  A canonical embedding  of  $\text{\bf  Aut}_{\text{\bf hol}}(\bf \mathbb C^n)$}
\author{{\sc Francisco Braun}\thanks{Partial support provided by the grant 2019/07316-0 of the S\~ao Paulo Research Foundation (FAPESP) and by a travel grant from Texas Christian University.}  and {\sc Frederico  Xavier}\thanks{Partial support provided by the John William and Helen Stubbs Potter Professorship.}}

\date{}

\maketitle

\begin{abstract}
{\small \noindent The group $\text{\rm  Aut}_{\text{hol}}(\mathbb C^n)$  of self-biholomorphisms of $\mathbb C^n$ consists of affine  maps if  $n=1$, but in higher dimensions it  is a  large   object that has not been described explicitly. Despite  the intricacies   involved when $n>1$,   surprisingly  every $F\in \text{\rm Aut}_{\text{hol}}(\mathbb C^n)$  is  uniquely   determined  inside the group by only two    data, of infinitesimal and global nature:  the $1$-jet  of $F$ at $0$, and the  complex Hessian    of  a certain  plurisubharmonic function  associated to $F$. If  $n=1$ this  global datum    is   zero for all  $F$, which is then  determined solely by its $1$-jet at $0$, and  one recovers  
$\text{\rm Aut}_{\text{hol}}(\mathbb C)= \text{\rm Aff}(\mathbb C)\cong \mathbb C  \times \mathbb C^{*}$.  Our main result,  formulated as the existence of a canonical embedding of $ \text{Aut}_{\text{hol}} ( \mathbb C^n)$, also   singles out    a natural   candidate  for   moduli space  of $ \text{Aut}_{\text{hol}} ( \mathbb C^n)$,     for all  $n>1$. }  
\end{abstract}

\section {Introduction}

Let  $\text{Aut}_{\text{hol}} ( \mathbb C^n)$ be  the group of all  holomorphic automorphisms  of $\mathbb C^n$,   and  $\text{\rm Aff} (\mathbb C^n)$  the subgroup of  invertible  affine  maps.  
It is elementary that   $\text{Aut}_{\text{hol}} ( \mathbb C)=\text{\rm Aff} (\mathbb C)\cong \mathbb C  \times \mathbb C^{*}$,  but in higher dimensions  $\text{Aut}_{\text{hol}} ( \mathbb C^n)$ is so large and has such a rich structure \cite {AL} that an explicit description of it is yet to be found.

In this paper we introduce  an  unified  approach to  the study of $ \text{Aut}_{\text{hol}} ( \mathbb C^n)$  that   puts   this  size dichotomy in context, and  at the same time identifies   in every  dimension  a natural   candidate  for  moduli space of $ \text{Aut}_{\text{hol}} ( \mathbb C^n)$.   

Our main result, Theorem \ref{proto-moduli},  establishes the existence of a canonical    embedding   
$$ \text{Aut}_{\text{hol}} ( \mathbb C^n) \stackrel{\frak{X}^n} {\longrightarrow}  \mathbb C^n \times GL(n, \mathbb C) \times \frak{L}^n(PSH), $$  
 where $\frak{L}^n(PSH)$ is the space of Levi matrices (complex Hessians) of real-analytic  plurisubharmonic functions on $\mathbb C^n-\{0\}$.
 
 The  embedding  $\frak{X}^n$ records   the $1$-jet  of the automorphism $F$ at $0$,  together with  a datum in $\frak{L}^n(PSH)$.  
When $n=1$ this global  datum is  zero for all  $F$,  and so all the information carried by $\frak{X}^n$ is already encoded in the $1$-jet of $F$ at $0$. From injectivity  one then recovers the elementary description
$\text{Aut}_{\text{hol}} ( \mathbb C)=\text{\rm Aff} (\mathbb C)$, with the corresponding moduli space $\mathbb C \times \mathbb C^{*}\times \text{\rm singleton}\cong \mathbb C \times \mathbb C^{*}$. 

The  range of  $\frak{X}^n$  can be regarded as a  ``proto-moduli space" \rm for the automorphism group, a concept that will be explained in the next section.  At present, what  stands between our  results and a full-fledged moduli space for $\text{Aut}_{\text{hol}} ( \mathbb C^n)$   is the lack of  an intrinsic characterization of  $\text{\rm Ran}(\frak{X}^n)$ when $n>1$.

Given the unusual nature of the map $\frak{X}^n$, it  merits commenting  on the process that led to its discovery.  In the course of  proving an embedding theorem for  $\text{Aut}_{\text{hol}} ( \mathbb C^n)$, at some point  of the argument one has to show that two holomorphic automorphisms $F$, $G$ are equal or, what is the same,  that $F\circ G^{-1}$ is the identity map. In this regard, the work \cite{X1} is relevant, as it contains  a characterization of  the identity  of $\mathbb C^n$  among injective local biholomorphisms. The \it statement \rm of Theorem \ref{proto-moduli} was arrived at - in a manner akin to reverse engineering -   by starting with  a  \it conjectural \rm  enhancement  of the main result of  \cite{X1}, and then working backwards to figure out what the adequate hypotheses  for Theorem \ref{proto-moduli} should be.  
  
Although it  may yet be possible to  establish   the desired  technical  improvement of \cite{X1},  in which case  Theorem \ref{proto-moduli} would be  a  corollary of  such a result,  it so happens that a sharper version of \cite{X1}   is no longer needed  to prove  our embedding theorem.  In an unexpected  turn of events, the present authors  realized that,  since \cite{X1} has  already led  to  the  right hypotheses in Theorem \ref{proto-moduli},   one can  actually  supply  an independent, non-technical, proof of this result.

The results in this paper are part of a  larger program  whose aim is to   use tools from differential geometry, analysis, topology,  and dynamical systems  in order to  understand the various  mechanisms behind the phenomenon of global injectivity of maps (for more on this,  see \cite {X2} and the references therein). 

\section{A proto-moduli space for $\text{Aut}_{\text{hol}} ( \bf \mathbb C^n)$}

In order to present  our findings about $\text{Aut}_{\text{hol}} ( \mathbb C^n)$
in a  systematic  way,  and  with a view towards future developments, 
 it is  convenient  to  adopt  a more formal approach to the   notion of ``moduli problem" for a  given ``space" $X$. This problem  is  generally  understood to  be  the search for a parametrization of $X$ by a set whose elements have  special, often geometric,  properties.  The following   set-theoretic   definitions make the idea of having ``special properties"  precise. Given  an injective map $f:X\to Y$ with range $Z$,  the bijection $Z\stackrel{f^{-1}}{\longrightarrow} X$     
  can be thought of as  a  \it  proto-parametrization  of $X$\rm,  and    $Z$ can be regarded as    a \it  proto-moduli space  of $X$\rm.

If, in  addition,   the membership   relation ``$y\in  Z$"  can  be described alternatively    using only   properties that  \it do not \rm involve  the set  $X$ (\it a fortiori\rm, invoking neither $f$ nor $Z$),   
 then  we  say that  \it $Z$ admits an intrinsic description\rm,     \it $Z\stackrel{f^{-1}}{\longrightarrow} X$     is a parametrization of $X$, and \it $Z$ is a moduli space of $X$\rm.    

To summarize,  the moduli spaces of $X$ are precisely  those ranges of embeddings of $X$ that admit an intrinsic description.  The desideratum that the moduli of $X$ should have ``special properties" is reflected in the condition that membership in $Z$ can be described  by properties that do not refer to $X$. A concrete example  illustrating  these simple  concepts  will be given later in this section, in connection with our main theorem.

The  advantage of  this 
 recasting   of  the   moduli problem    lies in the fact   that  the  problem   is   purposely  split into two  sub-problems, which can  then be examined separately:  one  first needs to   find  a  ``natural"  proto-parametrization  $f:X\to Y$  and then, in  a second stage,   argue   that the associated proto-moduli space    can be upgraded to  a full-fledged    moduli space.  
 
 In the present  context, the set $X$ is $\text{Aut}_{\text{hol}} ( \mathbb C^n)$.
Theorem  \ref{proto-moduli}  below provides  a solution of the first  half   of the moduli problem in all dimensions that, when specialized to $n=1$,    also gives  the expected solution of the second half  of the problem.  
We conjecture that these proto-moduli spaces for $\text{Aut}_{\text{hol}} ( \mathbb C^n)$  are  actual moduli spaces  when $n>1$ as well.

Before stating our results, we    review some basic concepts from the theory of several complex variables (a suitable reference is \cite[chapter 4]{H}). 
Given  a $C^2$ real-valued function $g$ defined  on  an open set of $\mathbb C^n$,   the  $n\times n$ Hermitean matrix 
${L}(g)=\Big (\frac{\partial^2 g}{\partial z_i \partial \overline{z}_j}  \Big )$
is  called the Levi matrix, or complex Hessian,  of $g$.  When $L(g)\geq 0$ (resp., $L(g)=0$) everywhere, $g$ is said to be plurisubharmonic (resp., pluriharmonic). 

Alternatively, the 
 plurisubharmonicity (pluriharmonicity) of $g$  can be characterized by  the property that   the restriction of $g$ to  every complex line intersecting the domain of $g$ is subharmonic (resp., harmonic).  Both notions  are invariant under holomorphic changes of coordinates.  A  prototypical example of a smooth plurisubharmonic function is  
 $\log \|G\|$, where $G:\mathbb C^m \to \mathbb C^n-\{0\}$ is  holomorphic. A pluriharmonic function defined on a simply connected set is the real part of a holomorphic function on this set.

We denote by 
 $\frak{L}^n(PSH)$ the  subset  of $C^{\omega}(\mathbb C^n-\{0\}, \mathbb C^{n\times n}) $
consisting of  the  Levi matrices of    real-analytic plurisubharmonic  functions defined on $\mathbb C^n-\{0\}$.

 \begin{thm}\label{proto-moduli}  
The  map $ {\frak{X}}^n:\text{\rm Aut}_{\text{\rm hol}} ( \mathbb C^n) \longrightarrow \mathbb C^n \times GL(n, \mathbb C) \times \frak{L}^n(PSH)$, 
 \begin{eqnarray}\label{embedding}
{ \frak{X}}^n(F)=\left(F(0), DF(0), L (\log \|DF(0)^{-1}(F-F(0))\|)\right), 
\end{eqnarray}
is injective for all $n\geq 1$. 
\end{thm} 

Let  $\frak{Q^n}$  denote the quotient space of the equivalence relation that identifies  two real-analytic plurisubharmonic  functions defined on $\mathbb C^n-\{0\}$ if their difference is pluriharmonic. There is a natural bijection between  $\frak{L}^n(PSH)$ and  $\frak{Q^n}$, given by $L(g)\to [g]$.  Accordingly, Theorem \ref{proto-moduli} can be stated  with a more intrinsic flavor:

\begin{thm} \label{proto-moduli(alt)} Any  $F\in \text{\rm Aut}_{\text{\rm hol}} ( \mathbb C^n)$ is uniquely determined inside this group by only two data: the $1$-jet of $F$ at $0$ and the equivalence class of $ \log \|DF(0)^{-1}(F-F(0))\|$ in $\frak{Q^n}$.
\end{thm}

The case $n=1$ is included in the statement of  Theorem \ref{proto-moduli}  for  the sake of uniformity,  but  it can be checked directly. Notice that  $\text{\rm Ran} (\frak{X}^1) \cong  \mathbb C \times \mathbb C^{*}\times \text {\rm singleton}$, a description that does not involve  $\text{\rm Aut}_{\text{\rm hol}} ( \mathbb C)$.
In particular,  $\text{\rm Ran}({\frak{X}}^n)$ is a  proto-moduli  space of $\text{\rm Aut}_{\text{\rm hol}} ( \mathbb C^n)$  for all $n> 1$,   and a moduli space   if  $n=1$ (see the discussion below for further details.)

\begin{conjecture} 
$\text{\rm Ran}(\frak{X}^n)$ is a  moduli  space  for all $n\geq 1$. 
\end{conjecture}

\noindent \rm The  central  issue here  is whether  $\text{\rm Ran} (\frak{X}^n)$   admits an  intrinsic  description   when  $n>1$, in which case 
we  would have a complete solution of the moduli problem for $\text{Aut}_{\text{hol}} ( \mathbb C^n)$, via the bijection 
$$\displaystyle  \text{\rm Ran} (\frak{X}^n) \stackrel{(\frak{X}^n)^{-1}}{\longrightarrow} 
\text{Aut}_{\text{hol}} ( \mathbb C^n).$$

For illustration purposes, let us  assume that $\frak{X}^1$ is injective, as stated in Theorem \ref{proto-moduli},     and  proceed to  show that $\text{\rm Ran} (\frak{X}^1)$  is a moduli space.  Along the way,  we will recover $\text{Aut}_{\text{hol}} ( \mathbb C)=\text{\rm Aff} (\mathbb C)$ and provide an informal reasoning that sheds light into the conceptual question of why $\text{Aut}_{\text{hol}} (\mathbb C^n)$ is so much bigger in higher dimensions.  As we shall see, ultimately,  this dichotomy is linked to a simple property of plurisubharmonic/pluriharmonic functions (Corollary \ref{affine}). 

Since  $\Delta \log |h|=0$ for any non-vanishing holomorphic function $h$ of a single variable, the  Levi matrix   in (\ref{embedding}) vanishes identically when $n=1$. It then follows  that $\frak{X}^1(F)=\frak{X}^1(F(0)+DF(0))$. 
As   $\frak{X}^1$  was assumed to be injective,   $F=F(0)+DF(0)$,   and  so     $\text{Aut}_{\text{hol}} ( \mathbb C)=\text{\rm Aff} (\mathbb C)$. 
Manifestly,  $\text{\rm Ran} (\frak{X}^1)\cong  \mathbb C \times \mathbb C^{*}\times \{0\}$,   a description independent of $\text{Aut}_{\text{hol}} ( \mathbb C)$, and so   $\mathbb C \times \mathbb C^{*}\times \{0\}$  qualifies as  a moduli space for $\text{Aut}_{\text{hol}} ( \mathbb C)$,  as per our    definition of moduli space.  The corresponding  parametrization of $\text{Aut}_{\text{hol}} ( \mathbb C)$ is 
$$\text{\rm Ran} (\frak{X}^1)\cong  \mathbb C \times \mathbb C^{*}\times \{0\} \stackrel{(\frak{X}^1)^{-1}} {\longrightarrow }\text{Aut}_{\text{hol}} ( \mathbb C), \;\;\;
(b, a, 0)\to \Phi_{(b,a,0)}, \;\;\; \Phi_{(b,a,0)}(z)=az+b.
$$

A  map $F \in \text{Aut}_{\text{hol}} ( \mathbb C^n)$ lies in  $\text{\rm Aff} (\mathbb C^n)$ if and only if 
$F=F(0)+DF(0)$. From Theorem \ref{proto-moduli} one has
\begin{cor}\label{affine} An automorphism $F\in \text{Aut}_{\text{hol}} ( \mathbb C^n)$ is affine if and only if  
\begin{eqnarray}\label{Eq} \log\left(\frac {\|DF(0)^{-1}(F-F(0))\|}{\| I \|}\right) \end{eqnarray} 
is pluriharmonic away from zero. 
\end{cor}

It is only when $n=1$  that condition (\ref{Eq}) is trivially satisfied, as the norm  becomes  absolute value, and 
\begin{eqnarray*}\label{Eq2} \Delta\log\left|\frac {DF(0)^{-1}(F(z)-F(0))}{ z }\right|=0, \;\; z\neq 0. \end{eqnarray*} 
Once again,  one recovers $\text{Aut}_{\text{hol}} ( \mathbb C)=\text{\rm Aff} (\mathbb C)$.  One can think of Corollary \ref{affine} as  providing an informal ``explanation" to why $ \text{Aut}_{\text{hol}} ( \mathbb C^n)$  is  much  larger than  $\text{\rm Aff} (\mathbb C^n)$   when  $n>1$.

 \section {Proof of the embedding theorem}
 
 We   write 
$\text{Aut}^0_{\text{hol}} ( \mathbb C^n)$ for the subgroup of $\text{Aut}_{\text{hol}} ( \mathbb C^n)$ consisting of  those $F$ that are   normalized at $0$ by   $F(0)=0$ and $DF(0)=I$. 
The map  $\Theta^n: \text{Aut}_{\text{hol}} ( \mathbb C^n) \to 
\text{Aut}^0_{\text{hol}} ( \mathbb C^n)$,      $$\Theta^n(F)=DF(0)^{-1}(F-F(0)),$$  is a retraction, in the sense that  the   restriction of $\Theta^n$ to $\text{Aut}^0_{\text{hol}} ( \mathbb C^n)$ is the identity map. 

Also, $\Theta^n(F)=\Theta^n(G)$ if and only if $F=H\circ G$ for some    $H\in \text{\rm Aff} (\mathbb C^n)$. In particular,   $\Theta^n$   induces a bijection  between a coset space and the subgroup of normalized  holomorphic automorphisms, namely, 
 $$\text{Aut}_{\text{hol}} ( \mathbb C^n)/\text{\rm Aff} (\mathbb C^n) \cong \text{Aut}^0_{\text{hol}} ( \mathbb C^n). $$ 
Thus, modulo affine maps, the study of 
$\text{Aut}_{\text{hol}} ( \mathbb C^n)$  is reduced to that of  $\text{Aut}^0_{\text{hol}} ( \mathbb C^n)$.

\begin{thm}\label{proto-moduli2}  
The map $ \displaystyle  \frak{X}_0^n:\text{\rm Aut}^0_{\text{\rm hol}} ( \mathbb C^n) \longrightarrow \frak{L}^n(PSH)$,  
$$\frak{X}_0^n(F)=L(\log \|F\|),$$ 
is injective for all $n\geq 1$. 
\end{thm}  
In words, Theorem \ref{proto-moduli2}  states that two normalized self-biholomorphisms  $F$ and $G$ of $\mathbb C^n$ are equal if and only if the difference between the plurisubharmonic functions $\log \| F \|$  and $\log \| G \|$ is pluriharmonic away from zero.

As in Theorem \ref{proto-moduli}, we include the case $n=1$ in the statement  of Theorem \ref{proto-moduli2} for  the sake of uniformity. In fact,  one  can check directly that both $\ \text{Aut}^0_{\text{hol}} ( \mathbb C) $ and 
$\text{Ran}  (\frak{X}_0^1)$ are singletons (i.e. the identity map of $\mathbb C $ versus   the zero $1\times 1$ matrix), and so the theorem holds trivially in this case.

Theorem \ref{proto-moduli} follows directly from Theorem  \ref{proto-moduli2}. Indeed, from
\begin{eqnarray} \label{frak-frak} \frak{X}^n(F)=\big(F(0), DF(0), \frak{X}_0^n(\Theta^n(F))\big)\end{eqnarray}
one  sees that  if  $\frak{X}^n_0$ is injective so is  $\frak{X}^n$ (and conversely).

To begin the proof of  Theorem  \ref{proto-moduli2}, assume 
$$\frak{X}_0^n(F)=L(\log \|F\|)=\frak{X}_0^n(G)=L(\log \|G\|).$$
Hence
$$L\Big(\log\frac{\|F\|}{\|G\|}\Big)(z)=0,  \;\;z\neq 0, $$
and so $\log\frac{\|F\|}{\|G\|}$ is pluriharmonic  for $z\neq 0$.
Because of the normalizations, $\frac{\|F\|}{\|G\|}$ extends continuously to $z=0$, with value $1$. If $n=1$ the   harmonic function $\log\frac{\|F\|}{\|G\|}$ is bounded on a punctured neighborhood of $0$, and so the singularity is removable.  In particular, $\log\frac{\|F\|}{\|G\|}$ is the real part of an entire function $f$. 

If $n>1$, the pluriharmonic function $\log\frac{\|F\|}{\|G\|}$ is defined on the simply-connected set $\mathbb C^n-\{0\}$, and so it  is the real part of a holomorphic function on $\mathbb C^n-\{0\}$. In higher dimensions  holomorphic functions have no isolated singularities, and so  one also obtains an entire function $f$ as before. 

Therefore, regardless of the dimension, one can write, 
$$\log\frac{\|F\|}{\|G\|}(z)=\text{Re} f(z),\;\;\;
\|F(z)\|=  |e^{f(z)}|\; \|G(z)\|,$$
with $f:\mathbb C^n \to \mathbb C$ holomorphic. 
Clearly,  $\text{Re} f(0)=0$  and by adding a suitable purely imaginary constant,  if necessary, we may assume that $f(0)=0$.  
Composing with $F^{-1}$ in the last equation, we have 
$$\|z\|=|e^{(f\circ F^{-1})(z)}|\; \|(G\circ F^{-1})(z)\|.$$
Setting  $H=e^{(f\circ F^{-1})}\; (G\circ F^{-1})$, one has $\|H(z)\|=\|z\|$.  
Hence,  the holomorphic map $H$ applies the unit ball into itself,   $H(0)=0$ and from $f(0)=0$ one has $DH(0)=I$.  
It now  follows from  a result  of Cartan \cite[p.66]{N} that  
$H=I$. Unwinding this condition  one has 
$$G\circ F^{-1}=hI, \;\; h=e^{-(f\circ F^{-1})}, \;\;h(0)=1.$$
In order to show that $F=G$,  it suffices to show that $h=1$ everywhere. 
It follows from  $G\circ F^{-1}=hI$ that  the map $hI$ is injective.  
Since the restriction of $hI$ to the line determined by $0$ and any  $z\in \mathbb C^n-\{0\}$ is injective, so is the map $\mathbb{C} \ni t \to h(tz) (t z)\in \mathbb{C}^n$. If $z_i$ is a non-zero coordinate of $z$,  a short argument  shows that  the function $\mathbb{C} \ni t \to h(tz) (t z_i)\in \mathbb{C}$  is injective as well. 
Since injective entire maps of $\mathbb{C}$ are affine linear, it follows that $h(tz)$ is constant as a function of $t$, and setting $t=0,1$ one has $h(z) = h(0)=1$. As $z\in \mathbb C^n-\{0\}$ was arbitrary, $h=1$ everywhere, thus concluding  the proof of Theorem \ref{proto-moduli2}.
\qed

In view of (\ref{frak-frak}),   the key to solving the conjecture stated in  section  2, hence the moduli problem for  $\text{\rm Aut}^0_{\text{\rm hol}} ( \mathbb C^n)$, is to find an intrinsic description of  $\text{Ran} ( \frak{X}_0^n)$ for $n>1$. 

\noindent 
Departamento de Matem\'{a}tica, Universidade Federal de S\~ao Carlos, S\~ao Carlos, SP, Brazil. 
\noindent franciscobraun@dm.ufscar.br

\vspace{.15cm}

\noindent Department of Mathematics, Texas Christian University, Fort Worth, TX, USA. 

\noindent f.j.xavier@tcu.edu

\end{document}